\documentclass[11pt, hidelinks]{article}
\usepackage[colorlinks=true,linkcolor=blue]{hyperref}
\usepackage{url}
\usepackage[margin=1in]{geometry}
\usepackage{amsmath,amssymb,amsthm, amsfonts}
\usepackage[T1]{fontenc}
\usepackage{mathtools}
\usepackage{xcolor}
\usepackage{enumitem, comment, xifthen}
\usepackage{graphicx}
\graphicspath{{figures/}}
\usepackage{subfig}

\theoremstyle{plain}

\newtheorem{theo}{Theorem}

\newtheorem{lem}{Lemma}
\newtheorem{prop}{Proposition}
\newtheorem{cor}{Corollary}

\theoremstyle{definition} 

\newtheorem{nota}{Notation}
\newtheorem{de}{Definition}
\newtheorem{exa}{Example}
\newtheorem{as}{Assumption}
\newtheorem{alg}{Algorithm}

\newcommand{\btheo}{\begin{theo}}
\newcommand{\bde}{\begin{de}}
\newcommand{\ble}{\begin{lem}}
\newcommand{\bpr}{\begin{prop}}
\newcommand{\bno}{\begin{nota}}
\newcommand{\bex}{\begin{exa}}
\newcommand{\bcor}{\begin{cor}}
\newcommand{\spro}{\begin{proof}}
\newcommand{\bas}{\begin{as}}
\newcommand{\balg}{\begin{alg}}

\newcommand{\etheo}{\end{theo}}
\newcommand{\ede}{\end{de}}
\newcommand{\ele}{\end{lem}}
\newcommand{\epr}{\end{prop}}
\newcommand{\eno}{\end{nota}}
\newcommand{\eex}{\end{exa}}
\newcommand{\ecor}{\end{cor}}
\newcommand{\fpro}{\end{proof}}
\newcommand{\eas}{\end{as}}
\newcommand{\ealg}{\end{alg}}

\theoremstyle{plain}

\newtheorem{theos}{Theorem}
\newtheorem{props}{Proposition}
\newtheorem{lems}{Lemma}
\newtheorem{cors}{Corollary}

\theoremstyle{definition}
\newtheorem{exas}{Example}
\newtheorem{algs}{Algorithm}
\newtheorem{asss}{Assumption}
\newtheorem{defns}{Definition}

\newcommand{\btheos}{\begin{theos}}
\newcommand{\etheos}{\end{theos}}
\newcommand{\bprops}{\begin{props}}
\newcommand{\eprops}{\end{props}}
\newcommand{\bdes}{\begin{defns}}
\newcommand{\edes}{\end{defns}}
\newcommand{\blems}{\begin{lems}}
\newcommand{\elems}{\end{lems}}
\newcommand{\bcors}{\begin{cors}}
\newcommand{\ecors}{\end{cors}}
\newcommand{\bexs}{\begin{exas}}
\newcommand{\eexs}{\end{exas}}
\newcommand{\balgs}{\begin{algs}}
\newcommand{\ealgs}{\end{algs}}
\newcommand{\bass}{\begin{asss}}
\newcommand{\eass}{\end{asss}}

\usepackage{hyperref}
\usepackage{url}
\usepackage[margin=1in]{geometry}
\usepackage{amsmath,amssymb,amsthm, amsfonts}
\usepackage[T1]{fontenc}
\usepackage{mathtools}
\usepackage{xcolor}
\usepackage{enumitem, comment, xifthen}
\usepackage{graphicx}
\usepackage{mathabx} 
\graphicspath{{figs/}}
\DeclarePairedDelimiter{\vecnorm}{\lVert}{\rVert}
\DeclarePairedDelimiter{\matnorm}{\vvvert}{\vvvert}
\newcommand{\twonorm}[2][]{\vecnorm[#1]{#2}_{2}}
\newcommand{\opnorm}[2][]{\matnorm[#1]{#2}_{\rm op}}

\makeatletter
\def\letterdef#1#2#3{\def\letterdef@##1{\expandafter\def\csname #1\endcsname{#2}}%
  \letterdef@@#3{?\@car{}}\@nil}
\def\letterdef@@#1{\@gobble#1\letterdef@{#1}\letterdef@@}
\makeatother

\DeclarePairedDelimiterX{\klx}[2]{(}{)}{%
  #1\;\delimsize\|\;#2%
}
\DeclarePairedDelimiterX{\quantklx}[3]{(}{)}{%
  #1\;\delimsize\|\;#2\;\delimsize\vert\;#3%
}
\DeclarePairedDelimiterX{\inner}[2]{\langle}{\rangle}{%
  #1,#2%
}

\newcommand{\R}{\mathbf R} 

\letterdef{c#1}{\mathcal{#1}}{ABCDEFGHIJKLMNOPQRSTUVWXYZ}
\letterdef{b#1}{\mathbb{#1}}{ABCDEFGHIJKLMNOPQRSTUVWXYZ}
\let\defn\coloneq
\newcommand{\twomax}[2]{\ensuremath{#1 \lor #2}}
\newcommand{\twomin}[2]{\ensuremath{#1 \land #2}}

\newcommand{\floor}[1]{\left\lfloor #1 \right\rfloor}

\newcommand{\e}{\mathrm{e}}

\newcommand{\ud}[0]{\,\mathrm{d}}  
\newcommand{\1}{\mathbf 1} 
\let\ones\1
 
\let\epsilon\varepsilon
\newcommand{\eps}{\varepsilon}

\let\tilde\widetilde

\let\succeq\succcurlyeq

\renewcommand{\leq}{\leqslant}
\renewcommand{\geq}{\geqslant}

\newcommand{\argmin}{\mathop{\rm arg\,min}}

\newcommand{\T}{\mathsf{T}}

\makeatletter
\newcommand{\E}{\operatorname*{\mathbf{E}}\ilimits@}
\makeatother
\makeatletter
\renewcommand{\P}{\operatorname*{\mathbf{P}}\ilimits@}
\makeatother

\newcommand{\ie}{\textit{i}.\textit{e}., }
\newcommand{\eg}{\textit{e}.\textit{g}., }

\newcounter{algorithmctr}
\renewcommand{\thealgorithmctr}{\arabic{algorithmctr}}
   {\refstepcounter{algorithmctr}\begin{list}{}{%
       \setlength{\rightmargin}{0\linewidth}%
       \setlength{\leftmargin}{0\linewidth}}%
       \rmfamily\small
       \item[]{\setlength{\parskip}{0ex}\hrulefill\par%
        \nopagebreak{\bfseries\textsf{Algorithm \thealgorithmctr~}}}}%
   {{\setlength{\parskip}{-1ex}\nopagebreak\par\hrulefill} \end{list}}

\makeatletter
\long\def\@makecaption#1#2{
        \vskip 0.8ex
        \setbox\@tempboxa\hbox{\small {\bf #1.} #2}
        \parindent 1.5em 
        \dimen0=\hsize
        \advance\dimen0 by -3em
        \ifdim \wd\@tempboxa >\dimen0
                \hbox to \hsize{
                        \parindent 0em
                        \hfil 
                        \parbox{\dimen0}{\def\baselinestretch{0.96}\small
                                {\bf #1.} #2
                                } 
                        \hfil}
        \else \hbox to \hsize{\hfil \box\@tempboxa \hfil}
        \fi
        }
\makeatother

\allowdisplaybreaks

\usepackage[colorlinks=true,linkcolor=blue]{hyperref}
\newcommand{\dimension}{d}
\newcommand{\numobs}{n}
\newcommand{\thetastar}{\theta^\star}

\newcommand{\Normal}[2]{\mathsf{N}\left(#1, #2\right)}
\newcommand{\ols}{{\widehat \theta^{\sf OLS}}}

\newcommand{\EmpCov}{\Sigma_\numobs}

\newcommand{\invsnr}{\tau} 
\newcommand{\effinvsnr}{\tau_\numobs}
\newcommand{\mustar}{\mu^\star}
\newcommand{\soft}[2]{\mathsf{S}_{#1}\!\!\left(#2\right)}
\newcommand{\Bin}[2]{\mathsf{Bin}\left(#1, #2\right)}
\newcommand{\Ber}[1]{\mathsf{Ber}\left(#1\right)}
\newcommand{\stols}[1]{{\widehat \theta^{\sf STOLS}_{#1}}}

\newcommand{\lambdathresh}{\overline{\lambda}}
\let\tilde\widetilde
\let\hat\widehat
\let\top\T

\newcommand{\upstairs}[1]{\textsuperscript{#1}}
\newcommand{\affilone}{\dag}

\newcommand{\affilchic}{$\diamond$}

\usepackage{color}
\definecolor{cm}{RGB}{0,0,200}
\definecolor{rp}{RGB}{0,200,0}

\begin{document}
\begin{center}
  {\bf{\Large On the design-dependent suboptimality of the Lasso}} \\
  \vspace*{.2in}
    \vspace*{.2in}

  \begin{tabular}{cc}
    Reese Pathak\upstairs{\affilone}\qquad \qquad Cong Ma\upstairs{\affilchic} \\[1.5ex]
    \upstairs{\affilone} Department of Electrical
    Engineering and Computer Sciences, UC Berkeley \\
    \upstairs{\affilchic} Department of Statistics, University of
    Chicago \\ 
    \texttt{pathakr@berkeley.edu, congm@uchicago.edu}
  \end{tabular}
  \vspace*{.2in}
\end{center}

\begin{abstract}
This paper investigates the effect of the design matrix on the ability (or inability) to estimate a sparse parameter in linear regression. 
More specifically, we characterize the optimal rate of estimation when the smallest singular value of the design matrix is bounded away from zero.
In addition to this information-theoretic result, we provide and analyze a procedure which is simultaneously statistically optimal and computationally efficient, based on soft thresholding the ordinary least squares estimator. 
Most surprisingly, we show that the Lasso estimator---despite its widespread adoption for sparse linear regression---is provably minimax rate-suboptimal 
when the minimum singular value is small. 
We present a family of design matrices and sparse parameters for which we can guarantee that the Lasso with \emph{any} choice of regularization parameter---including those which are data-dependent and randomized---would fail in the sense that 
its estimation rate is suboptimal by polynomial factors in the sample size. 
Our lower bound is strong enough to preclude the statistical optimality of all forms of the Lasso, including its highly popular penalized, norm-constrained, and cross-validated variants. 
\end{abstract}

\section{Introduction}\label{sec:intro}

In this paper, we consider the standard linear regression model 
\begin{equation}\label{eq:linear-regression}
y = X \thetastar + w, 
\end{equation}
where $\thetastar \in \R^{\dimension}$ is the unknown parameter, $X \in \R^{\numobs \times \dimension}$ is the design matrix, and $w \sim \Normal{0}{\sigma^2I_\numobs}$ denotes the stochastic noise. Such linear regression models are pervasive in statistical analysis~\cite{lehmann2006theory}. 
To improve model selection and estimation, it is often desirable to impose a \emph{sparsity} 
assumption on $\thetastar$---for instance, 
we might assume that $\thetastar$ has few nonzero entries or that 
it has few large entries. This amounts to assuming that for some $p \in [0,1]$
\begin{equation}\label{eq:sparsity-assumption}
	 \|\thetastar\|_p \leq R,
\end{equation}
where $\|\cdot \|_p$  denotes the $\ell_{p}$ vector (quasi)norm, and $R > 0$ is the radius of the $\ell_p$ ball. 
There has been a flurry of research on this sparse 
linear regression model~\eqref{eq:linear-regression}-\eqref{eq:sparsity-assumption} over the last three decades; see the recent books~\cite{buhlmann2011statistics, hastie2015statistical, wainwright2019high,fan2020statistical,Johnstone2019} for an overview. 

Comparatively less studied, is the effect of the design matrix $X$ on the 
ability (or inability) to estimate $\thetastar$ under the sparsity assumption. 
Intuitively, when $X$ is ``close to singular'', 
we would expect that certain directions of $\thetastar$ would be 
difficult to estimate. Therefore, in this paper we seek to determine 
the optimal rate of estimation when the \emph{smallest singular value of $X$ is bounded.}  
More precisely, we consider the following set of design matrices 
\begin{equation}\label{eqn:covariate-shift-X-class}
\cX_{\numobs, \dimension}(B) \defn \Big\{\, X \in \R^{\numobs \times
  \dimension} : \frac{1}{\numobs} X^\T X \succeq \frac{1}{B} I_{d}\, \Big\}, 
\end{equation}
and aim to characterize the corresponding minimax rate of estimation
\begin{equation}\label{eq:minimax-risk}
\mathfrak{M}_{\numobs, d}(p, \sigma, R, B)
\defn
\inf_{\hat \theta} \sup_{\substack{X \in \cX_{\numobs,\dimension}(B) \\ 
 \|\thetastar\|_p \leq R}}
\E_{y \sim \Normal{X\thetastar}{\sigma^2 I_\numobs}} \Big[\vecnorm[\big]{\hat \theta - \thetastar}^2_2 \Big]. 
\end{equation}

\subsection{A motivation from learning under covariate shift}

Although it may seem a bit technical to focus on the dependence of the estimation error on the 
smallest singular value of the design matrix $X$, we would like to point out an additional motivation 
which is more practical and also motivates our problem formulation. 
This is the problem of linear regression in a well-specified model with covariate shift. 

To begin with, recall that under random design, in the standard linear observational model (\ie without covariate shift) the statistician observes random covariate-label pairs of the form $(x,y)$. Here, the covariate $x$ is drawn from a distribution $Q$ and the label $y$ satisfies 
$\E[y \mid x] = x^\top \thetastar$. 
The goal is to find an estimator $\hat \theta$ that minimizes the out-of-sample excess risk, which takes the quadratic form $\E_{x \sim Q}[ ((\hat \theta - \thetastar)^\T x)^2]$. 
When the covariate distribution $Q$ is isotropic, meaning that $\E_{x \sim Q} [x x^\top] = I$, 
the out-of-sample excess risk equals the squared $\ell_2$ error 
$\|\hat \theta - \thetastar\|_2^2$.

Under covariate shift, there is a slight twist to the standard linear regression model previously described, where now the covariates $x$ are drawn from a (source) distribution $P$ that 
differs from the (target) distribution $Q$ under which we would like to deploy our estimator. 
Assuming $Q$ is isotropic, the goal is therefore still to minimize the out-of-sample excess risk under $Q$, which is $\|\thetastar - \hat \theta\|_2^2$. 
In general, if $P \neq Q$ and no additional assumptions are made, then learning with covariate shift is impossible in the sense that no estimator can be consistent for the 
optimal parameter $\thetastar$. 
It is therefore common (and necessary) to impose some additional assumptions on the pair $(P, Q)$ 
to facilitate learning. One popular assumption relates to the likelihood ratio between the source-target pair. It is common to assume that absolute continuity holds so that $Q \ll P$ and 
that the the likelihood ratio $\tfrac{\ud Q}{\ud P}$ is uniformly bounded~\cite{ma2023optimally}. 
Interestingly, it is possible to show that if $\tfrac{\ud Q}{\ud P}(x) \leq B$ for 
$P$-almost every $x$, then the semidefinite inequality \begin{equation}\label{ineq:cov-shift}
\E_{x \sim P}[ x x^\top] \succeq \frac{1}{B} I
\end{equation}
holds~\cite{ma2023optimally, wang2023pseudo}.
Comparing the inequality~\eqref{ineq:cov-shift} to our class $\cX_{\numobs, \dimension}(B)$ 
as defined in display~\eqref{eqn:covariate-shift-X-class}, 
we note that our setup can be regarded as a fixed-design variant of linear regression with covariate shift~\cite{lei2021near, eyre2023out, zhang2022class}.

\subsection{Determining the minimax rate of estimation}
\label{sec:minimax-rate}

We begin with one of our main results, which precisely characterizes the (order-wise) minimax risk $\mathfrak{M}_{\numobs, d}(p, \sigma, R, B)$ of estimating $\thetastar$ under the sparsity constraint $\|\thetastar\|_p \leq R$ and over the restricted design class $\cX_{\numobs, \dimension}(B)$.

\btheo \label{thm:minimax-rate}
Let $n \geq d \geq 1$ and $\sigma, R, B > 0$ be given, and put 
$\effinvsnr^2 \defn \tfrac{\sigma^2 B}{R^2 \numobs}$.
There exist two universal constants 
$c_\ell, c_u$ satisfying $0 < c_\ell < c_u < \infty$ such that
\begin{enumerate}[label=(\alph*)]
\item \label{thm:rate-weakly-sparse}
if $p \in (0, 1]$ and $\effinvsnr^2 \in [d^{-2/p}, \log^{-1}(\e d)]$, then 
\[
c_\ell \, R^2\,  \Big(\effinvsnr^2 \log\big(\e d \effinvsnr^p\big)\Big)^{1 - p/2}
\leq 
\mathfrak{M}_{\numobs, d}(p, \sigma, R, B)
\leq 
c_u \, R^2 \, \Big(\effinvsnr^2 \log\big(\e d \effinvsnr^p\big)\Big)^{1 - p/2}
, 
\quad \mbox{and}
\]
\item \label{thm:rate-hard-sparse}
if $p = 0$, we denote $s = R \in [\dimension]$, and have 
\[
c_\ell \, \frac{\sigma^2 B}{\numobs} s \log\Big(\e \frac{d}{s}\Big)
\leq 
\mathfrak{M}_{\numobs, d}(p, \sigma, s, B)
\leq 
c_u \, \frac{\sigma^2 B}{\numobs} s \log\Big(\e \frac{d}{s}\Big).
\]
\end{enumerate}
\etheo 

\medskip
\noindent The proof of Theorem~\ref{thm:minimax-rate}  relies on a reduction to the Gaussian sequence model~\cite{Johnstone2019}, and is deferred to 
Section~\ref{sec:proof-thm-minimax-rate}.

Several remarks on Theorem~\ref{thm:minimax-rate} are in order. 
The first observation is that Theorem~\ref{thm:minimax-rate}
 is sharp, apart from universal constants that do not depend on the tuple of problem parameters~$(p, n, d, \sigma, s, R, B)$. 

Secondly, it is worth commenting on the sample size restrictions in Theorem~\ref{thm:minimax-rate}. For all $p \in [0, 1]$, we have assumed the ``low-dimensional'' setup that the number of observations $n$ dominates the dimension $d$.\footnote{Notably, this still allows 
$n$ to be proportional to $d$, \eg we can tolerate $n = d$.}
Note that this is necessary for the class of designs 
$\cX_{\numobs, \dimension}(B)$ to be nonempty. 
On the other hand, for $p > 0$ we additionally 
require that the sample size is ``moderate'', i.e., $\effinvsnr^2 \in [d^{-2/p}, \log^{-1}(\e d)]$. 
We make this assumption so that we can focus on what we believe is 
the ``interesting'' regime: where neither 
ordinary least squares nor constant estimators are optimal. 
Indeed, when $n \geq d$ but $\effinvsnr^2 \geq 
\log^{-1}(\e d)$, it is easily verified that the optimal rate of 
estimation is on the order $R^2$; intuitively the effective noise level 
is too high and no estimator can dominate $\hat \theta \equiv 0$ uniformly.  
On the other hand, when $n \geq d$ but $\effinvsnr^2 \leq d^{-2/p}$, then 
the ordinary least squares estimator is minimax optimal; intuitively, the noise level is sufficiently small such that there is, in the worst case, no need to shrink on the basis of the $\ell_p$ constraint to achieve the optimal rate.

Last but not least, as shown in Theorem~\ref{thm:minimax-rate}, the optimal rate of estimation depends on the signal-to-noise ratio $\effinvsnr^{-2} = \numobs R^2 / (\sigma^2 B)$. 
As $B$ increases, the design $X$ becomes closer to singular, estimation of $\thetastar$, as expected, 
becomes more challenging.  
The dependence of our result on $B$ is exactly analagous to the impact of the likelihood ratio bound $B$ appearing in the context of prior 
work on nonparametric regression under covariate shift~\cite{ma2023optimally}.

%


\subsection{A computationally efficient estimator}
\label{sec:STOLS}

The optimal estimator underlying the proof of Theorem~\ref{thm:minimax-rate} 
requires computing a $\dimension$-dimensional Gaussian integral, and 
therefore is not computationally efficient in general. 
In this section we propose an estimator that is both  
computationally efficient and statistically optimal, up to constant factors. 

Our procedure is based on the soft thresholding operator: for $v \in \R^\dimension$ and $\eta > 0$, 
we define 
\[
\soft{\eta}{v} \defn \argmin_{u \in \R^\dimension}
\Big\{\,
\|u - v\|_2^2 + 2 \eta \|u\|_1
\,\Big\}.
\]
Note that soft thresholding involves a coordinate-separable optimization problem
and has an explicit representation, thus allowing efficient computation.
Then we define the \emph{soft thresholded ordinary least squares estimator} 
\begin{equation}\label{eqn:soft-ols}
\stols{\eta}(X, y) \defn 
\soft{\eta}{\ols(X, y)},
\end{equation}
where $\ols(X, y)$ is the usual ordinary least squares estimate---equal to 
$(X^\T X)^{-1} X^\T y$ in our case. 
We have the following guarantees for its performance. 
\btheo \label{prop:performance-of-soft-thresholded OLS}
The soft thresholded ordinary least squares estimator~\eqref{eqn:soft-ols} satisfies 
\begin{enumerate}[label=(\alph*)]
\item \label{thm-part:soft-ols-weak-sparse}
in the case $p \in (0, 1]$, for any $R > 0$, 
if $\effinvsnr^2 \in [d^{-2/p}, \log^{-1}(\e d)]$, then 
\[
\sup_{X \in \cX_{\numobs, \dimension}(B)}
\sup_{\|\thetastar\|_p \leq R}
\E\Big[
\|\stols{\eta}(X, y) - \thetastar\|_2^2\Big]  
\leq 
6 \, R^2 \, 
\Big(\effinvsnr^2 \log(\e d \effinvsnr^p)\Big)^{1 -p/2},
\]
with the choice $\eta = \sqrt{2 R^2 \effinvsnr^2 \log(\e d \effinvsnr^p)}$, and 
\item \label{thm-part:soft-ols-hard-sparse}
in the case $p = 0$, for any $s \in [\dimension]$, 
\[
\sup_{X \in \cX_{\numobs, \dimension}(B)}
\sup_{\|\thetastar\|_0 \leq s}
\E\Big[
\|\stols{\eta}(X, y) - \thetastar\|_2^2\Big] 
\leq 
6 \, \frac{\sigma^2 B}{\numobs} s \log\Big(\e \frac{d}{s}\Big),
\]
with the choice $\eta = \sqrt{2 \frac{\sigma^2 B}{\numobs} \log(\tfrac{\e \dimension}{s})}$.
\end{enumerate}
\etheo 

\medskip
\noindent The proof is presented in Section~\ref{sec:proof-stols}.

Comparing the guarantee in Theorem~\ref{prop:performance-of-soft-thresholded OLS} to the minimax rate in 
Theorem~\ref{thm:minimax-rate}, it is immediate to see that the soft thresholded ordinary 
least squares estimator is  minimax optimal apart from constant factors. 

Secondly, we would like to point out a (simple) modification to the soft thresholding ordinary least squares procedure that allows it to be adaptive to the hardness of the 
particular design matrix encountered. 
To achieve this, note that $X \in \cX_{\numobs, \dimension}(\hat B)$ for $\hat B \defn \opnorm{(X^\T X)^{-1}}$. Therefore the results in Theorem~\ref{prop:performance-of-soft-thresholded OLS} continue to hold 
with $B$ replaced by (a possibly smaller) $\hat B$, provided that the thresholding parameter~$\eta$ is properly adjusted. For instance, 
in the case with $p = 0$, we have 
\begin{equation}\label{ineq:sharper-bound-stols-hard-sparse}
\sup_{\|\thetastar\|_0 \leq s}
\E\Big[
\|\stols{\hat \eta}(X, y) - \thetastar\|_2^2\Big] 
\leq 
6 \, \frac{\sigma^2 \hat B}{\numobs} s \log\Big(\e \frac{d}{s}\Big),
\end{equation}
provided we take 
$\hat \eta = 
\sqrt{2 \frac{\sigma^2 \hat B}{\numobs} \log(\tfrac{\e \dimension}{s})}$.

Finally, we note that inspecting our proof, the upper bound for $\stols{\hat \eta}(X, y)$ also holds for a larger set of design matrices
\begin{equation*}
\cX^{\mathsf{diag}}_{\numobs, \dimension}(B) \defn \Big\{\, X \in \R^{\numobs \times
  \dimension} : \left( \frac{1}{\numobs} X^\T X \right)^{-1}_{ii} \leq {B}, \quad \text{for }1\leq i \leq d\, \Big\}.
\end{equation*}
Since $\cX_{\numobs, \dimension}(B) \subset \cX^{\mathsf{diag}}_{\numobs, \dimension}(B)$, this means 
after combining the lower bounds in Theorem~\ref{thm:minimax-rate} with the guarantees in 
Theorem~\ref{prop:performance-of-soft-thresholded OLS}, we additionally have established the minimax rate 
over this larger family $\cX^{\mathsf{diag}}_{\numobs, \dimension}(B)$.
%
%

\subsection{Is Lasso optimal?}

Arguably, the Lasso estimator~\cite{tibshirani1996regression} is the most widely used estimator for sparse linear regression. 
Given a regularization parameter $\lambda > 0$, the Lasso is defined to be
\begin{equation}\label{eqn:def-lasso}
\hat \theta_\lambda(X, y) \defn \argmin_{\vartheta \in \R^\dimension}
\Big\{\, \frac{1}{\numobs} 
\twonorm{X\vartheta - y}^2 + 2 \lambda \vecnorm{\vartheta}_1\,\Big\}.
\end{equation}
Surprisingly, we show that the Lasso estimator---despite its popularity---is provably \emph{suboptimal} for estimating $\thetastar$ when $B \gg 1$. 
\bcor\label{cor:lasso-example}
The Lasso is minimax suboptimal by polynomial factors in the sample size when $d = n$ and 
$B = \sqrt{\numobs}$.
More precisely, 
\begin{enumerate}[label=(\alph*)] 
\item if $p \in (0, 1]$, and $\sigma = R = 1$, 
then we have 
\[
\sup_{X \in \cX_{\numobs,\dimension}(B)}
\sup_{\|\thetastar\|_p \leq R}
\E\Big[\inf_{\lambda > 0} \|\hat \theta_\lambda(X, y) - \thetastar\|_2^2 
\Big]
\gtrsim 
1, \quad \text{and}
\]
\item if $p = 0$, and $\sigma = s = 1$, then we have
\[
\sup_{X \in \cX_{\numobs,\dimension}(B)}
\sup_{\|\thetastar\|_0 \leq s}
\E\Big[\inf_{\lambda > 0} \|\hat \theta_\lambda(X, y) - \thetastar\|_2^2 
\Big]
\gtrsim 1. 
\]
\end{enumerate}
\ecor

\medskip 
\noindent Corollary~\ref{cor:lasso-example} is in fact a special case of a more general theorem (Theorem~\ref{thm:suboptimality-penalized-lasso}) to be provided later. 

Applying Theorem~\ref{thm:minimax-rate} to the regime considered in Corollary~\ref{cor:lasso-example},
 we obtain the optimal rate of estimation 
 \[
\Big(\frac{\sqrt{1 + \log \numobs}}{\numobs^{1/4}}\Big)^{2 - p}
\ll 1, \qquad \text{for every } p \in [0,1].
\]
As shown, in the worst-case, the multiplicative gap between the 
performance of the Lasso and a minimax optimal estimator in this 
scaling regime is at least polynomial in the sample size.
As a result, the Lasso is quite strikingly
minimax \emph{suboptimal} in this scaling regime.

In fact, the lower bound against Lasso in Corollary~\ref{cor:lasso-example} is extremely strong. 
Note that in the lower bound, the Lasso is even allowed to leverage the oracle information $\thetastar$ to 
calculate the optimal instance-dependent choice of the 
regularization parameter (\emph{c.f.,} $\inf_{\lambda > 0} \|\hat \theta_\lambda(X, y) - \thetastar\|_2^2$). 
As a result, the lower bound applies to  
any estimator which can be written as the penalized Lasso estimator 
with data-dependent choice of penalty.
Many typical Lasso-based estimators, such as the norm-constrained 
and cross-validated Lasso, can be written as the penalized Lasso 
with a data-dependent choice of the penalty parameter $\lambda$. 
For instance, in the case of the norm-constrained Lasso, this holds by convex duality. 
Thus, we can rule out the minimax optimality 
of any procedure of this type, in light of Corollary~\ref{cor:lasso-example}. 

The separation between the oracle Lasso and the minimax optimal estimator 
can also be demonstrated in experiments, as shown below in 
Figure~\ref{fig:simulation}. 

\begin{figure}[h!t]
\centering 
\hspace{-9mm}
\includegraphics[width=0.35\linewidth]{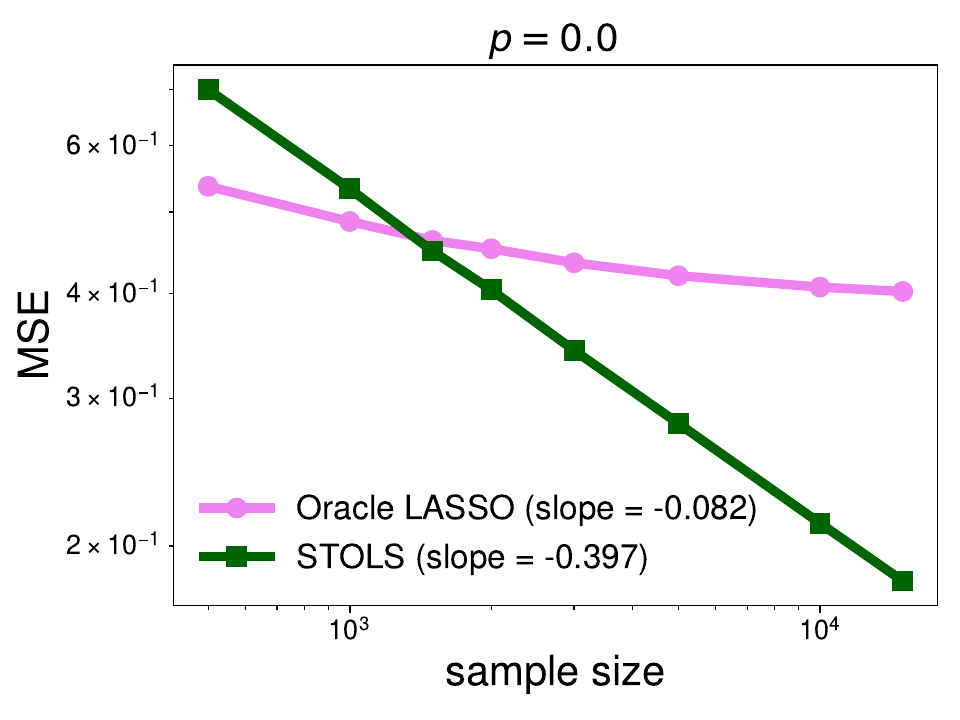}
\hspace{-3mm}
\includegraphics[width=0.35\linewidth]{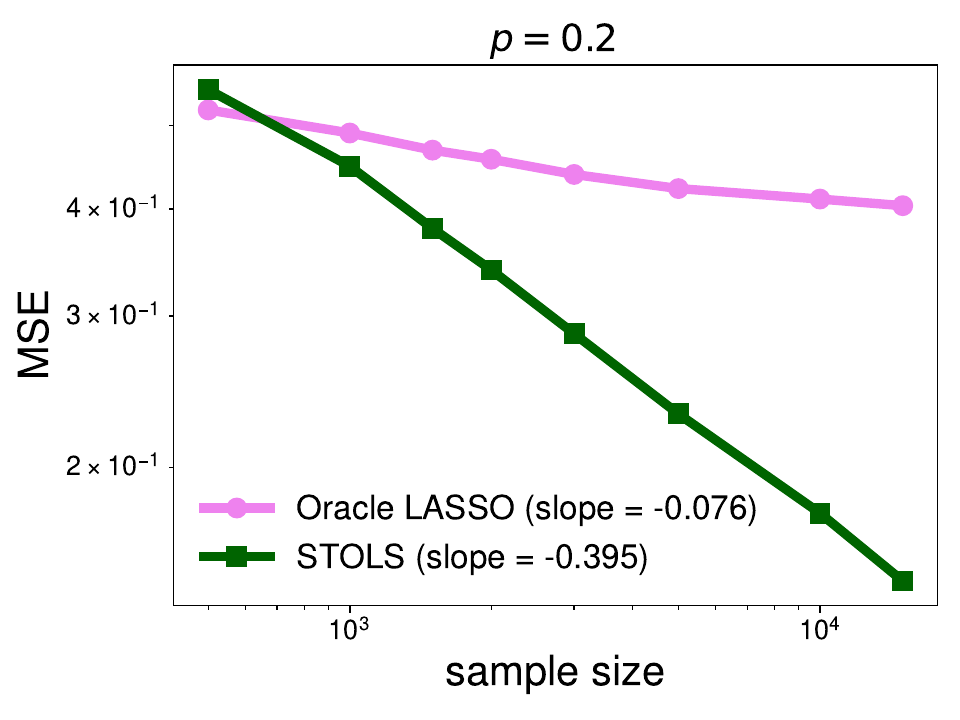}
\hspace{-3mm}
\includegraphics[width=0.35\linewidth]{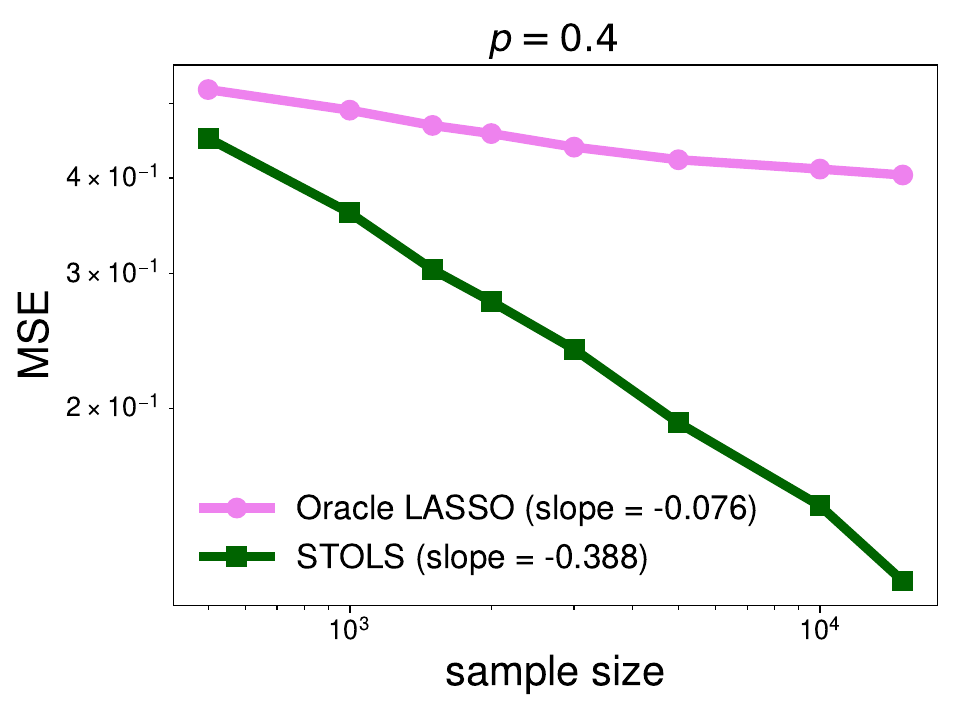}
\caption{Numerical simulation demonstrating suboptimality of the Lasso, 
with the oracle choice of regularization $\hat \lambda 
\in \argmin_{\lambda > 0} \|\hat \theta_\lambda(X, y) - \thetastar\|_2^2$ 
versus the soft thresholded ordinary least squares (STOLS) procedure as defined in
display~\eqref{eqn:soft-ols}. We have simulated the lower bound instance from Corollary~\ref{cor:lasso-example}; 
for each sample size $\numobs$, 
we simulate oracle Lasso and STOLS on a pair $(X_\numobs, \thetastar_\numobs)$, with the dimension $d = \numobs$ and lower singular 
value bound $B = \sqrt{\numobs}$. Further details on the simulation are provided in Appendix~\ref{app:sim-details}.}
\label{fig:simulation}
\end{figure}

\subsection{Connections to prior work}

In this section, we draw connections and comparisons between our work and existing literature. 
\paragraph{Linear regression with elliptical or no constraints.} 
Without any parameter restrictions, the exact minimax rate for linear regression
when error is measured in the $\ell_2$ norm along with the dependence on the design matrix is known: it is 
given by $\sigma^2 \mathsf{Tr}( (X^\top X)^{-1})$~\cite{lehmann2006theory}. 
These results match our intuition that as the smallest singular value of $X$ decreases, the hardness of estimating $\thetastar$ increases. 
It is also worth mentioning that the design matrix does not play a role, apart from being invertible, in determining the optimal rate for the in-sample prediction error. The rate is given uniformly by $\frac{\sigma^2 \dimension}{ \numobs}$ when $n \geq d$~\cite{hsu2012random}.   

On the other hand, with $\ell_2$- or elliptical parameter  constraints the minimax rate in both fixed and random design 
was established in the recent paper~\cite{Pat23}. 
Although that work shows the dependence on the design, 
the rate is not explicit and the achievable result requires potentially solving a semidefinite program. 
More explicit results were previously derived under proportional asymptotics in the paper~\cite{Dic16}, 
in the restricted setting of Gaussian isotropic random design. 
The author is able to establish the asymptotic minimaxity of a particular ridge regression estimator.
These type of results are not immediately useful for our work, since they are based on linear shrinkage 
procedures which are known to be minimax-suboptimal even in orthogonal design, in the $\ell_p$ setting 
for $p < 2$~\cite{DonJoh94}.

\paragraph{Gaussian sequence model and sparse linear regression.} In the case of orthogonal design, \ie when $\frac{1}{\numobs} X^\T X =  I_{d}$, the minimax risk of estimating sparse $\thetastar$ is known; see~\cite{DonJoh94}. It can also be shown that Lasso, with optimally-tuned parameter, can achieve the (order-wise) optimal estimation rate~\cite{Johnstone2019}. 
This roughly corresponds to the case with $B= 1$ in our consideration. Our work generalizes this line of research in the sense that we characterize the optimal rate of estimation over the larger family of design matrices $\frac{1}{\numobs} X^\T X \succeq \frac{1}{B} I_{d}$. In stark contrast to the Gaussian sequence model, Lasso is no longer optimal, even with the oracle knowledge of the true parameter.  

Without assuming an orthogonal design, \cite{candes2013well} provides a design-dependent lower bound in the exact sparsity case (\ie $p=0$). The lower bound depends on the design through its Frobenius norm $\|X\|_{\mathsf{F}}^2$.  Similarly, in the weak sparsity case, \cite{raskutti2011minimax} provides lower bounds depending on the maximum column norm of the design matrix.  However, matching upper bounds are not provided in this general design case. In contrast, using the minimum singular value of $X$ (\emph{c.f.,} the parameter $B$) allows us to obtain matching upper and lower bounds in sparse regression.

\paragraph{Suboptimality of Lasso.} The suboptimality of Lasso for 
minimizing the \emph{prediction} error has been noted in the case of
exact sparsity (\ie $p=0$). 
To our knowledge, previous studies required a carefully chosen 
design matrix which was highly-correlated. For instance, 
it was shown that for certain highly-correlated designs the Lasso can 
achieve only a slow rate $(1/\sqrt{n})$, while information-theoretically, the optimal rate is faster ($1/n$); see for instance the papers~\cite{van2018tight, Kelner2021preconditioning}. Additionally in the paper~\cite{candes2009near}, the authors exhibit the failure of Lasso for a \emph{fixed} regularization parameter, which does not necessarily rule out the optimality of 
other Lasso variants. Similarly, in the paper~\cite{foygel2011fast}, it is shown via a correlated design matrix and a 2-sparse vector, that the  norm-constrained version of Lasso can only achieve a slow rate in terms of the prediction error. Again, this result does not rule out the optimality of 
other variants of the Lasso. In addition, in the paper~\cite{DALALYAN17Lasso}, there is an example for which Lasso with any fixed (\ie independent from the observed data) choice of regularization would fail to achieve the optimal rate. Again, this fails to rule out data-dependent choices of regularization or other variants of the Lasso. In our work, we are able to rule out the optimality of the Lasso by considering a simple diagonal 
design matrix which exhibits no correlations among the columns. 
Nonetheless, for any $p \in [0, 1]$, we show that the Lasso will fall short of optimality by polynomial factors in the sample size. Our result 
also simultaneously rules out the optimality of constrained, penalized, and 
even data-dependent variants of the Lasso, in contrast to the literature described above. 

\paragraph{Covariate shift.}
As mentioned previously, our work is also related to linear regression under covariate shift~\cite{lei2021near, zhang2022class, eyre2023out}. The statistical analysis of covariate shift, albeit with an asymptotic nature, dates back to the seminal work by Shimodaira~\cite{shimodaira2000improving}.  Recently, nonasymptotic minimax analysis of covariate shift has gained much attention in unconstrained parametric models~\cite{ge2023maximum}, nonparametric classification~\cite{kpotufe2021marginal}, and also nonparametric regression~\cite{pathak2022new, ma2023optimally, wang2023pseudo}.

%
%
%

%
%
%
%
%
%

\section{A closer look at the failure mode of Lasso}
In this section, we take a closer look at the failure instance for Lasso. 
We will investigate the performance of the Lasso on diagonal 
design matrices $X_\alpha \in \R^{\numobs \times \dimension}$
which satisfy, when $d = 2k$, 
\[
\frac{1}{\numobs} X_\alpha^\T X_\alpha = 
\begin{pmatrix} 
\tfrac{\alpha}{B} I_k & 0 \\ 
0 & \tfrac{1}{B} I_k
\end{pmatrix}.
\]
Thus, this matrix has condition number $\alpha$ and satisfies 
$X_\alpha \in \cX_{\numobs, \dimension}(B)$ for all $\alpha \geq 1$.
As our proof of Theorem~\ref{thm:minimax-rate} reveals, 
from an information-theoretic perspective, the hardest design matrix
$X_\alpha$ is with the choice $\alpha = 1$: when all directions 
have the worst possible signal-to-noise ratio. Strikingly, 
this is not the case for the Lasso: 
there are in fact choices of $\alpha \gg 1$ which 
are even harder for the Lasso.

\btheo
\label{thm:suboptimality-penalized-lasso}
Fix $\numobs \geq \dimension \geq 2$ and let $\sigma, B > 0$ be given.
For $\alpha \geq 1$, on the diagonal design $X_\alpha$, 
\begin{enumerate}[label=(\alph*)]
\item if $p \in (0, 1]$ and $R > 0$, then there is a vector 
$\thetastar \in \R^\dimension$ such that $\|\thetastar\|_p \leq R$ but 
\[
\E_{y \sim \Normal{X_\alpha \thetastar}{\sigma^2 I_\numobs}}
\Big[\inf_{\lambda > 0} \|\hat \theta_\lambda(X_\alpha, y) - \thetastar\|_2^2\Big]
\geq 
\frac{9}{20000}\,
\Big(
\frac{\sigma^2 B \dimension}{\numobs \alpha} 
\wedge 
R^2\Big(\frac{\sigma^2 B}{R^2 \numobs} \alpha\Big)^{1 - p/2}
\wedge 
R^2
\Big), \quad \mbox{and}
\]
\item if $p = 0$ and $s \in [d]$, then there is a vector 
$\thetastar \in \R^d$ which is $s$-sparse but 
\[
\E_{y \sim \Normal{X_\alpha \thetastar}{\sigma^2 I_\numobs}}
\Big[\inf_{\lambda > 0} \|\hat \theta_\lambda(X_\alpha, y) - \thetastar\|_2^2\Big]
\geq 
\frac{9}{20000}\,
\Big(
\frac{\sigma^2 B \dimension}{\numobs \alpha} 
\wedge \frac{\sigma^2 B s}{\numobs} \alpha\Big).
\]
\end{enumerate}
\etheo

\medskip 
\noindent 
The proof of Theorem~\ref{thm:suboptimality-penalized-lasso} is presented
in Section~\ref{sec:proof-theorem-suboptimality-penalized-lasso}. 
We now make several comments on the implications of this result.

We emphasize that the dependence of the Lasso on 
the parameter $\alpha$, which governs
the condition number of the matrix $X_\alpha$, is suboptimal, 
as revealed by Theorem~\ref{thm:suboptimality-penalized-lasso}. 
At a high-level, large $\alpha$ should only make the ability to 
estimate $\thetastar$ \emph{easier}---it effectively increases the 
signal-to-noise ratio in certain directions. This can also 
be seen from Theorem~\ref{thm:minimax-rate}: the conditioning of the 
design matrix \emph{does not} enter into the worst-case rate of estimation 
when the bottom signular value of $X$ is bounded. Nonetheless, Theorem~\ref{thm:suboptimality-penalized-lasso} shows that the Lasso actually can 
suffer when the condition number $\alpha$ is large.

\paragraph{Proof of Corollary~\ref{cor:lasso-example}.}
We now complete the proof of Corollary~\ref{cor:lasso-example} given Theorem~\ref{thm:suboptimality-penalized-lasso}.

Maximizing over the parameter $\alpha \geq 1$ appearing in our 
result, we can determine a particularly nasty configuration of the 
conditioning of the design matrix for the Lasso. 
Doing so, we find that for $p \in (0, 1]$ and $R > 0$ that 
\[
\sup_{X \in \cX_{\numobs,\dimension}(B)}
\sup_{\|\thetastar\|_p \leq R}
\E\Big[\inf_{\lambda > 0} \|\hat \theta_\lambda(X, y) - \thetastar\|_2^2 
\Big]
\gtrsim 
R^2 \bigg(\twomin{\Big(\frac{\sigma^2 B \sqrt{d}}{R^2 \numobs}\Big)^{\tfrac{4 - 2p}{4 - p}}}{1}\bigg).
\]
This is exhibited by considering the lower bound in Theorem~\ref{thm:suboptimality-penalized-lasso} with the choice 
$\alpha^\star(p) = (\effinvsnr^2 d^{2/p})^{p/(4-p)}$.
On the other hand, if $p = 0$, we have for $s \in [d]$ that
\[
\sup_{X \in \cX_{\numobs,\dimension}(B)}
\sup_{\|\thetastar\|_0 \leq s}
\E\Big[\inf_{\lambda > 0} \|\hat \theta_\lambda(X, y) - \thetastar\|_2^2 
\Big]
\gtrsim 
\frac{\sigma^2 B}{\numobs} \sqrt{sd}
\]
The righthand side above is exhibited by considering the lower bound 
with the choice $\alpha^\star(0) = \sqrt{d/s}$.

The proof is completed by setting $d = n$, 
$B = \sqrt{\numobs}$, and $\sigma = R = 1$.

\section{Proofs}
In this section, we present the proofs for the main results of this paper. 
We start with introducing a few useful notations.
For a positive integer $k$, we define
$[k] \defn \{1, \dots, k\}$. For a real number $x$, we define 
$\floor{x}$ to be the largest integer less than or equal to $x$ and 
$\{x\}$ to be the fractional part of $x$.

\subsection{Proof of Theorem~\ref{thm:minimax-rate}}
\label{sec:proof-thm-minimax-rate}

Our proof is based on a decision-theoretic reduction to the Gaussian 
sequence model.
It holds in far greater generality, and so we actually prove 
a more general claim which could be of interest to other linear 
regression problems on other parameter spaces or with other loss functions. 

To develop the claim, we first need to introduce notation. Let 
$\Theta \subset \R^\dimension$ denote a parameter space, 
and let $\ell \colon \Theta \times \R^\dimension \to \R$ be a given loss 
function. We define two minimax rates, 
\begin{subequations}\label{eqn:defn-general-minimax-rates}
\begin{align}
\mathfrak{M}_{\rm seq}(\Theta, \ell, \nu) 
&\defn 
\inf_{\hat \mu} 
\sup_{\mustar \in \Theta} 
\E_{y \sim \Normal{\mustar}{\nu^2 I_\dimension}}
\Big[\ell\big(\mustar, \hat \mu(y)\big)\Big], \quad \mbox{and} 
\label{defn:minimax-rate-seq}
\\ 
\mathfrak{M}_{\rm reg}(\Theta, \ell, \sigma, B, \numobs)
&\defn 
\inf_{\hat \theta}
\sup_{X \in \cX_{\numobs, \dimension}(B)}
\sup_{\thetastar \in \Theta}
\E_{y \sim \Normal{X\thetastar}{\sigma^2 I_\numobs}}
\E\Big[\ell\big(\thetastar, \hat\theta(X, y)\big)\Big].
\label{defn:minimax-rate-reg}
\end{align}
\end{subequations}

The definitions above correspond to the minimax rates of estimation 
over the $\ell_p$ ball of radius $R > 0$ in $\R^\dimension$
for the Gaussian sequence model, in the case of definition~\eqref{defn:minimax-rate-seq}, and for $n$-sample 
linear regression with $B$-bounded design, 
in the case of definition~\eqref{defn:minimax-rate-reg}. The infima 
range over measurable functions of the observation vector $y$, in both cases. 

The main result we need is the following statistical reduction from 
linear regression to mean estimation in the Gaussian sequence model. 
\bpr [Reduction to sequence model] \label{prop:reduction}
Fix $n, d \geq 1$ and $\sigma, B > 0$. Let 
$\Theta \subset \R^\dimension$, and 
$\ell \colon \Theta \times \R^\dimension \to \R$ be given. 
If $\ell(\theta, \cdot) \colon \R^\dimension \to \R$ is a 
convex function for each $\theta \in \Theta$, then 
\[
\mathfrak{M}_{\rm reg}(\Theta, \ell, \sigma, B, \numobs)
=
\mathfrak{M}_{\rm seq}\bigg(\Theta, \ell, \sqrt{\frac{\sigma^2 B}{\numobs}}\bigg).
\]
\epr 

Deferring the proof of Proposition~\ref{prop:reduction} to 
Section~\ref{sec:proof-of-prop-reduction} for the moment, we note that it immediately implies Theorem~\ref{thm:minimax-rate}. 
Indeed, we set $\Theta = \Theta_{d, p}(R)$ and $\ell = \ell_{\rm sq}$ where 
\[
\Theta_{d, p}(R) \defn 
\{\theta \in \R^\dimension : \|\theta\|_p\leq R\}, \quad \mbox{and} \quad 
\ell_{\rm sq}(\theta, \hat\theta) = 
\|\hat \theta - \theta\|_2^2.
\]
With these choices, we obtain
\[
\mathfrak{M}_{\numobs, \dimension}(p, \sigma, R, B)
= 
\mathfrak{M}_{\rm reg}(\Theta_{d, p}(R), \ell_{\rm sq}, \sigma, B, \numobs)
= 
\mathfrak{M}_{\rm seq}\bigg(\Theta_{d, p}, \ell_{\rm sq}, 
\sqrt{\frac{\sigma^2 B}{\numobs}}\bigg),
\]
where the final equality follows from Proposition~\ref{prop:reduction}. 
The righthand side then corresponds to estimation in the $\ell_2$ norm over the Guassian sequence model with parameter space corresponding to an 
$\ell_p$ ball in $\R^\dimension$, which is a well-studied problem~\cite{DonJoh94,BirMas01,Johnstone2019}; 
we thus immediately obtain the result via classical results (for a precise statement with explicit constants, see Propositions~\ref{prop:minimax-rate-weak-sparse-gsm} and~\ref{prop:minimax-rate-hard-sparse-gsm}
presented in Appendix~\ref{app:gaussian-sequence-model}).

\subsubsection{Proof of Proposition~\ref{prop:reduction}}
\label{sec:proof-of-prop-reduction}

We begin by lower bounding the regression minimax risk by the sequence model 
minimax risk. Indeed, let $X_\star$ be such that $\tfrac{1}{\numobs} X_\star^\T X_\star = \tfrac{1}{B} I_\dimension$. Then we have 
\begin{align*}
\mathfrak{M}_{\rm reg}(\Theta, \ell, \sigma, B, \numobs)
&\geq 
\inf_{\hat \theta} \sup_{\thetastar \in \Theta} 
\E_{y \sim \Normal{X_\star \thetastar}{\sigma^2 I_\numobs}} 
\Big[ \ell\big(\thetastar, \hat \theta \big)\Big] \\ 
&\geq 
\inf_{\hat \theta} \sup_{\thetastar \in \Theta} 
\E_{z \sim \Normal{\thetastar}{\tfrac{\sigma^2 B}{\numobs} I_\dimension}} 
\Big[ \ell\big(\thetastar, \hat \theta \big)\Big] \\
&= \mathfrak{M}_{\rm seq}\bigg(\Theta, \ell_{\rm sq}, 
\sqrt{\frac{\sigma^2 B}{\numobs}}\bigg).
\end{align*}
The penultimate equality follows by noting that in the regression model, 
$\cP_{X_\star} \defn 
\{ \Normal{X_\star \theta}{\sigma^2 I_\dimension} : \theta \in \Theta\}$, 
the ordinary least squares (OLS) estimate is a sufficient statistic. 
Therefore, by the Rao-Blackwell Theorem, there exists a minimax optimal 
estimator which is only a function of the OLS estimate. 
For any $\thetastar$, the ordinary least squares estimator has the distribution, $\Normal{\thetastar}{\frac{\sigma^2 B}{\numobs} I_\dimension}$, which 
provides this equality. 

We now turn to the upper bound. Let $\ols(X, y) = (X^\T X)^{-1}X^\T y$ 
denote the ordinary least squares estimate. For any estimator $\hat \mu$, 
we define 
\begin{equation}\label{eqn:upper-bound-reduction-estimator}
\tilde \theta(X, y) 
\defn 
\E_{\xi \sim \Normal{0}{W}} \Big[
\hat \mu\big(\ols(X, y) + \xi \big)\Big], 
\qquad \mbox{where}\quad  W = \frac{\sigma^2 B}{\numobs} I_\dimension
- \sigma^2 (X^\T X)^{-1}.
\end{equation}
Note that for any $X \in \cX_{\numobs, \dimension}(B)$, we have 
$W \succeq 0$. 
Additionally, by Jensen's inequality, 
for any $X \in \cX_{\numobs, \dimension}(B)$, and any $\thetastar \in \Theta$, 
\begin{align*}
\E_{y \sim \Normal{X\thetastar}{\sigma^2 I_\numobs}} 
\Big[\ell\big(\thetastar, \tilde \theta(X, y)\big)\Big] 
&\leq 
\E_{y \sim \Normal{X\thetastar}{\sigma^2 I_\numobs}}
\E_{\xi \sim \Normal{0}{W}}
\Big[\ell\big(\thetastar, \hat \mu\big(\ols(X, y) + \xi \big)\big)\Big] \\
&= 
\E_{z \sim \Normal{\thetastar}{\tfrac{\sigma^2 B}{\numobs} I_\dimension}} 
\Big[\ell\big(\thetastar, \hat \mu(z)\big)\Big]
\end{align*}
Passing to the supremum over $\thetastar \in \Theta$ on each side and 
then taking the infimum over measurable estimators, we immediately see that
the above display implies 
\[
\mathfrak{M}_{\rm reg}(\Theta, \ell, \sigma, B, \numobs)
\leq \mathfrak{M}_{\rm seq}\bigg(\Theta, \ell, 
\sqrt{\frac{\sigma^2 B}{\numobs}}\bigg),
\]
as needed.

\subsection{Proof of Theorem~\ref{prop:performance-of-soft-thresholded OLS}}
\label{sec:proof-stols}
We begin by bounding the risk for soft thresholding procedures, 
based on a rescaling and monotonicity argument 
and applying results from~\cite{Johnstone2019}.
To state it, we need to define 
the quantitites 
\[
(\nu_\numobs^\star)^2 \defn \frac{\sigma^2 B}{\numobs}
\quad \mbox{and} 
\quad 
\rho_\numobs(\theta, \eta)
\defn 
\dimension \, (\nu^\star_\numobs)^2 
\e^{-(\eta/\nu^\star_\numobs)^2/2}
+ 
\sum_{i=1}^\dimension 
\twomin{\theta_i^2}{(\nu_\numobs^\star)^2} 
+ 
\sum_{i=1}^\dimension
\twomin{\theta_i^2}{\eta^2}.
\]
Then we have the following risk bound.

\ble \label{lem:stols-uniform-X-bound}
For any $\thetastar\in \R^\dimension$ and any $\eta > 0$ we have 
\[
\sup_{X \in \cX_{\numobs, \dimension}(B)} 
\E\Big[\twonorm{\stols{\eta}(X, y) - \thetastar}^2\Big]
\leq 
\rho_\numobs(\thetastar, \eta).
\]
\ele 

We now define for $\zeta > 0$ and a subset $\Theta \subset \R^\dimension$, 
\[
T(\zeta, \Theta) \defn 
\sup_{\thetastar \in \Theta}
\sum_{i=1}^\dimension 
\twomin{(\thetastar_i)^2}{\zeta^2}.
\]
Lemma~\ref{lem:stols-uniform-X-bound} then yields
with the choice $\eta = \gamma \nu_\numobs^\star$ for some 
$\gamma \geq 1$ that 
\begin{equation}\label{ineq:uniform-X-and-theta-bound-stols}
\sup_{\thetastar \in \Theta}
\sup_{X \in \cX_{\numobs, \dimension}(B)} 
\E\Big[\twonorm{\stols{\eta}(X, y) - \thetastar}^2\Big]
\leq 
3 \bigg[\twomax{d (\nu_\numobs^\star)^2 \e^{-\gamma^2/2}}{T(\gamma \nu^\star_\numobs, \Theta)}\bigg].
\end{equation}

We bound the map $T$ for the $\ell_p$ balls of interest. 
To state the bound, we 
use the shorthand $\Theta_p$ for the radius-$R$ $\ell_p$ ball in 
$\R^\dimension$ centered at the origin for $p \neq 0$, and 
for $p = 0$, the set of $s$-sparse vectors in $\R^\dimension$, 
for $s \in [\dimension]$. 
 
\ble
\label{lem:uniform-param-bounds-upper-soft}
Let $d \geq 1$ be fixed. We have the following relations:
\begin{enumerate}[label=(\alph*)]
  \item \label{ineq:weak-sparse-term-bounds-upper-soft}
  in the case $p \in (0, 1]$, we have for any $\zeta > 0$,
  \[
  T(\zeta, \Theta_p)
  \leq 
   R^2\,
  \bigg[ \Big(\frac{\zeta}{R}\Big)^{2} \dimension  
   \wedge \Big(\frac{\zeta}{R}\Big)^{2-p} \wedge 1 \bigg]
  \]
  for any $R > 0$, and 
  \item \label{ineq:hard-sparse-term-bounds-upper-soft}
  in the case $p = 0$, we have for any $\zeta > 0$,
  \[
  T(\zeta, \Theta_p)
  = \zeta^2 s,
  \]
  for any $s \in [\dimension]$.
\end{enumerate}
\ele 

To complete the argument, we now split into the two cases of hard and 
weak sparsity. 
\paragraph{When $p = 0$:} Combining inequality~\eqref{ineq:uniform-X-and-theta-bound-stols} together with Lemma~\ref{lem:uniform-param-bounds-upper-soft}, we find for 
$\eta = \gamma \nu_\numobs^\star$, $\gamma \geq 1$, that 
\[
\sup_{\thetastar \in \Theta}
\sup_{X \in \cX_{\numobs, \dimension}(B)} 
\E\Big[\twonorm{\stols{\eta}(X, y) - \thetastar}^2\Big]
\leq 
3 (\nu_\numobs^\star)^2 \Big[  \twomax{d \e^{-\gamma^2/2}}{\gamma^2 s} \Big]
= 6 (\nu_\numobs^\star)^2 s 
\log\Big(\e \frac{d}{s}\Big),
\]
where the last equality holds with $\gamma^2 = 2 \log(\e d/s)$. 

\paragraph{When $p \in (0, 1]:$} 
Combining inequality~\eqref{ineq:uniform-X-and-theta-bound-stols} together with Lemma~\ref{lem:uniform-param-bounds-upper-soft}, we find for
$\eta = \gamma \nu_\numobs^\star$, $\gamma \geq 1$ 
\begin{equation}\label{ineq:risk-bound-gamma-weak-sparse}
\sup_{\thetastar \in \Theta}
\sup_{X \in \cX_{\numobs, \dimension}(B)} 
\E\Big[\twonorm{\stols{\eta}(X, y) - \thetastar}^2\Big]
\leq 
3 R^2\, \bigg[
\twomax{d \effinvsnr^2 \e^{-\gamma^2/2}}{
  \Big(
  \gamma^{2 -p} \effinvsnr^{2-p} 
  \wedge 
  1\Big)
}\bigg]
\end{equation}
Above, we used $\effinvsnr^2 R^2 = (\nu_\numobs^\star)^2$ and 
$\gamma^2 \effinvsnr^2 \dimension \geq \gamma^{2-p} \effinvsnr^{2-p}$, 
which holds since $\gamma \geq 1$ and $\effinvsnr^2 \geq d^{-2/p}$.
If we take $\gamma^2 = 2 \log(\e \dimension \effinvsnr^{p})$, then 
note $\gamma^2 \geq 1$ by $\effinvsnr^2 \geq d^{-2/p}$ and
 the term in brackets in inequality~\eqref{ineq:risk-bound-gamma-weak-sparse} satisfies 
\begin{align*}
\twomax{d \effinvsnr^2 \e^{-\gamma^2/2}}{
  \Big(
  \gamma^{2 -p} \effinvsnr^{2-p} 
  \wedge 
  1\Big)}
  &= 
  \twomax{\frac{\effinvsnr^{2-p}}{\e}}{
   \Big( 
   (2 \effinvsnr^{2} \log(\e \dimension \effinvsnr^{p}))^{1-p/2} 
  \wedge 
  1\Big)} \\ 
  &\leq 
  2\bigg[\twomax{\effinvsnr^{2-p}}{
  \Big(  
   (\effinvsnr^{2} \log(\e \dimension \effinvsnr^{p}))^{1-p/2} 
  \wedge 
  1\Big)}\bigg]\\
  &=2  (\effinvsnr^{2} \log(\e \dimension \effinvsnr^{p}))^{1-p/2},
\end{align*}
which follows by $\effinvsnr^2 \in [d^{-2/p}, \log^{-1}(\e d)]$.

Thus, to complete the proof of Theorem~\ref{prop:performance-of-soft-thresholded OLS} we only need to provide the 
proofs of the lemmas used above.

\subsubsection{Proof of Lemma~\ref{lem:stols-uniform-X-bound}}

Note that if $z = \ols(X, y)$ then 
$\stols{\eta}(X, y) = \soft{\eta}{z} = \soft{\eta}{\thetastar + \xi}$ where 
$\xi \sim \Normal{0}{\tfrac{\sigma^2}{\numobs} \EmpCov^{-1}}$, 
where we recall $\EmpCov \defn (1/\numobs) X^\T X$. 
We now recall some classical results regarding the soft thresholding 
estimator. Let us write for $\lambda > 0$ and $\mu \in \R$,
\begin{align*}
r_{\rm S}(\lambda, \mu) &\defn 
\E_{y \sim \Normal{\mu}{1}} \Big[\big(\soft{\lambda}{y} - \mu\big)^2\Big], 
\quad \mbox{and,} \\
\tilde{r}_{\rm S}(\lambda, \mu) 
&\defn 
\e^{-\lambda^2/2} + 
\Big(\twomin{1}{\mu^2}\Big)
+ 
\Big(\twomin{\lambda^2}{\mu^2}\Big).
\end{align*}
Using $\twomin{(a+b)}{c} \leq \twomin{a}{c} + \twomin{b}{c}$ for nonnegative $a, b, c \geq 0$, Lemma 8.3 and the inequalities
$r_{\rm S}(\lambda, 0) \leq 1 + \lambda^2$ and 
$r_{\rm S}(\lambda, 0) \leq \e^{-\lambda^2/2}$ on page 219 of 
the monograph~\cite{Johnstone2019}, we find that 
$r_{\rm S}(\lambda, \mu)
\leq 
\tilde{r}_{\rm S}(\lambda, \mu).$
Define $\nu_i^2 \defn \frac{\sigma^2}{\numobs} (\EmpCov^{-1})_{ii}$. 
Using the fact that 
\(
\left(\soft{\eta}{z}\right)_i = \soft{\tfrac{\eta}{\nu_i}}{\frac{z_i}{\nu_i}}
\) 
for $i \in [\dimension]$, we obtain 
\[
\E[(\left(\soft{\eta}{z}\right)_i - \thetastar_i)^2]
= \nu_i^2 r_{\rm S}\Big(\frac{\eta}{\nu_i}, \frac{\thetastar_i}{\nu_i}\Big)
\leq 
\nu_i^2 \tilde{r}_{\rm S}\Big(\frac{\eta}{\nu_i}, \frac{\thetastar_i}{\nu_i}\Big).
\]
Summing over the coordinates yields 
\begin{align*}
\E\Big[\twonorm{\stols{\eta}(X, y) - \thetastar}^2\Big] 
&\leq \sum_{i=1}^\dimension
\nu_i^2 \tilde{r}_{\rm S}\Big(\frac{\eta}{\nu_i}, \frac{\thetastar_i}{\nu_i}\Big) \\ 
&= 
\sum_{i=1}^\dimension 
\nu_i^2 \e^{-(\eta/\nu_i)^2/2} 
+ \big(\twomin{(\thetastar_i)^2}{\nu_i^2}\big) 
+ \big(\twomin{(\thetastar_i)^2}{\eta^2}\big) \\ 
&\leq \rho_\numobs(\thetastar, \eta),
\end{align*}
where the last inequality follows by noting that 
both $\nu \mapsto \nu^2 \e^{-(\eta/\nu)^2/2}$ and 
$\nu \mapsto \twomin{\theta^2}{\nu^2}$ are nondecreasing 
functions of $\nu > 0$. Noting that this inequality holds uniformly on $X \in \cX_{\numobs, \dimension}(B)$ and passing to the supremum yields the 
claim.

\subsubsection{Proof of Lemma~\ref{lem:uniform-param-bounds-upper-soft}}

The proof of claim~\ref{ineq:hard-sparse-term-bounds-upper-soft} is 
immediate, so we focus on the case $p \in (0, 1]$, $R > 0$.
We consider three cases for the tuple $(R, \zeta, p, d)$. 
Combination of all three cases will yield the claim. 

\paragraph{When $R \geq \zeta d^{1/p}$:} 
Evidently, for each $\theta$ such that $\|\theta\|_p \leq R$, we have
\[
\sum_{i=1}^\dimension \twomin{\theta_i^2}{\zeta^2}
\leq
\zeta^2 \dimension.
\]

\paragraph{When $R \leq \zeta$:}
This case is immediate, since $\theta \in \Theta_p$ implies
$\|\theta\|_2 \leq \|\theta\|_p \leq R \leq \zeta$. 

\paragraph{When $\zeta \leq R \leq \zeta d^{1/p}$:}
In this case, by rescaling and putting $\eps \defn \tfrac{\zeta^2}{R^2}$, 
we have
\[
\sup_{\|\theta\|_p \leq R} 
\sum_{i=1}^\dimension \twomin{\theta_i^2}{\zeta^2}
=
R^2 \bigg[
\sup_{\lambda \in \Delta_\dimension} 
\sum_{i : \lambda_i \geq \eps^{p/2}} 
\eps 
+ 
\sum_{i : \lambda_i < \eps^{p/2}} 
\lambda_i^{2/p}\bigg]
= 
\eps R^2\Big(\floor{\eps^{-p/2}} + \{\eps^{-p/2}\}^{2/p}\Big)
\]
where above $\Delta_\dimension$ denotes the probability 
simplex in $\R^\dimension$. Noting that $\eps \leq 1$ and
$p \leq 1$ we have 
\[
\Big(\floor{\eps^{-p/2}} + \{\eps^{-p/2}\}^{2/p}\Big)
\leq 
\eps^{-p/2},
\]
which in combination with the previous display shows that
\[
\sup_{\|\theta\|_p \leq R} 
\sum_{i=1}^\dimension \twomin{\theta_i^2}{\zeta^2}
\leq
R^2 
\eps^{1-p/2}.
\]
To conclude, now note that $R^2 \eps^{1-p/2} = 
R^p \zeta^{2-p}$.

%

\subsection{Proof of Theorem~\ref{thm:suboptimality-penalized-lasso}}
\label{sec:proof-theorem-suboptimality-penalized-lasso}

Since $X_\alpha$ has nonzero entries only on the diagonal, we can derive 
an explicit representation of the Lasso estimate, as defined in display~\eqref{eqn:def-lasso}. To develop this, we first recall the notion of the 
\emph{soft thresholding operator}, which is defined by a 
parameter $\eta > 0$ and then satisfies
\[
\soft{\eta}{v} 
\defn 
\argmin_{u \in \R} \Big\{\, 
(u - v)^2 + 2 \eta |u|
\,\Big\}. 
\] 
We then start by stating the following lemma which is crucial for our analysis. It is a straightforward consequence of 
the observation that 
\begin{equation}\label{eqn:lasso-for-lower-bound-instance}
\hat \theta_\lambda(X_\alpha, y)_i = 
\begin{cases} 
\soft{\lambda B/\alpha}{z_i} 
& 1 \leq i \leq k \\ 
\soft{\lambda B}{z_i} 
& k + 1 \leq i \leq \dimension
\end{cases},
\end{equation}
where we have defined the independent random variables
$z_i \sim \Normal{\thetastar_i}{\frac{\sigma^2 B}{n \alpha}}$ if $i \leq k$ and
$z_i \sim \Normal{\thetastar_i}{\frac{\sigma^2 B }{n}}$ otherwise.
\ble 
Let $\thetastar \in \R^\dimension$. Then for the design matrix $X_\alpha$,
we have 
\[
\|\hat \theta_\lambda(X_\alpha, y) - \thetastar\|_2^2 
=
\sum_{i=1}^k
\Big(\soft{\lambda B/\alpha}{z_i} - \thetastar_i\Big)^2
+ 
\sum_{i=k+1}^\dimension
\Big(\soft{\lambda B}{z_i} - \thetastar_i\Big)^2. 
\]
\ele 

We will now focus on vectors $\thetastar(\eta) = (0_k, \eta, 0_{d - 2k})$, 
which are parameterized by $\eta \in \R^k$.
For these vectors, we can further lower bound the best risk as 
\begin{equation}\label{ineq:lower-terms-both}
\inf_{\lambda > 0} 
\|\hat \theta_\lambda(X_\alpha, y) - \thetastar(\eta)\|_2^2 
\geq 
\twomin{T_1}{T_2(\eta)} 
\end{equation}
where we have defined 
\[
\lambdathresh = \sqrt{\frac{\sigma^2}{\numobs}\frac{\alpha}{B}}, 
\quad  
T_1 \defn \inf_{\lambda \leq \lambdathresh} 
\sum_{i=1}^k
\Big(\soft{\lambda B/\alpha}{z_i}\Big)^2 
\quad \mbox{and} 
\quad 
T_2(\eta) \defn 
\inf_{\lambda \geq \lambdathresh}
\sum_{i=k+1}^{2k}
\Big(\soft{\lambda B}{z_i} - \eta_i\Big)^2. 
\]

We now move to lower bound $T_1$ and $T_2(\eta)$ by auxiliary, independent 
random variables.
\ble [Lower bound on $T_1$]
\label{lem:lower-bound-one}
Then, for any $\eta \in \R^k$ if $\thetastar = \thetastar(\eta)$, we have 
\[
T_1 \geq \frac{1}{4} \frac{\sigma^2 B}{\numobs \alpha}
 Z \qquad 
 \mbox{where} 
 \quad 
 Z \defn 
 \sum_{i=1}^k 
 \1\Big\{\tfrac{|z_i|}{\sqrt{\sigma^2 B/(n\alpha)}} \geq 3/2 \Big\}.
\]
\ele 

\ble [Lower bound on $T_2(\eta)$]
\label{lem:lower-bound-two}
Fix $\eta \in \R^k$ if $\thetastar = \thetastar(\eta)$, and suppose that 
\[
0 \leq \eta_i \leq 2 \sqrt{\frac{\sigma^2 B \alpha}{n}} \quad 
\mbox{for all}~i \in [k].
\] 
Then, we have 
\[
T_2(\eta) \geq \frac{1}{4} 
\sum_{i=1}^k \eta_i^2 W_i 
\quad 
\mbox{where} \quad 
W_i \defn \1\{z_{k + i} \leq \eta_i\}
\]
\ele 

Note that $Z$ is distributed as a Binomial random variable: 
$Z \sim \Bin{k}{p}$ where $p \defn \P\{|\Normal{0}{1}\!| \geq 3/2\}$. 
Similarly, $W_i$ are Bernoulli: we have $W_i \sim \Ber{1/2}$.

\paragraph{Lower bound for $p > 0$:} 
We consider two choices of $\eta$. First suppose that 
$4 \effinvsnr^2 \alpha \leq 1$. 
Then, we will consider $\eta = R \delta \1_\ell$ where 
\[
\delta \defn 2 \effinvsnr \sqrt{\alpha}
\quad \mbox{and} \quad 
\ell \defn \twomin{k}{\floor{\delta^{-p}}}
\]
For this choice of $\eta$ we have by assumption that 
$\effinvsnr^2 \geq d^{-2/p}$ that $\ell \geq (1/2) \delta^{-p}$
and so 
\[
T_2(\eta) \geq \frac{1}{4} R^2 \delta^2 \ell  \overline{W_\ell} 
\geq \frac{R^2}{8} \delta^{2 - p}  \overline{W_\ell} 
\geq \frac{R^2}{8} 
 \Big(\effinvsnr^2 \alpha \Big)^{1 - p/2} \overline{W_\ell} 
\]
Above, $\overline{W_\ell} \defn (1/\ell)\sum_{i=1}^\ell W_i$.
On the other hand, if $4 \effinvsnr^2 \alpha \geq 1$, we take 
$\eta' = R e_1$, and we consequently obtain 
\[
T_2(\eta') \geq \frac{R^2}{4} W_1
\]
Taking $\delta = 1/2$ in Lemma~\ref{lem:lower-bound-one}, let us define 
\[
c_1 \defn \min_{1 \leq \ell \leq k} 
\P\Big\{\overline{W_\ell} \geq \frac{1}{2}\Big\} 
\quad \mbox{and} \quad 
c_2 \defn \P\Big\{ Z \geq k p\Big\}.
\]
Let us take 
\[
\thetastar_\alpha \defn \begin{cases}
\thetastar(\eta) & 4 \effinvsnr^2 \alpha \leq 1 \\ 
\thetastar(\eta') & 4 \effinvsnr^2 \alpha > 1 
\end{cases}.
\]
Then combining Lemmas~\ref{lem:lower-bound-one} and~\ref{lem:lower-bound-two}
and the lower bounds on $T_2(\eta), T_2(\eta')$ above, we see that 
\begin{equation}\label{ineq:final-lower-weak-sparse-lasso}
\E_{y \sim \Normal{X_\alpha \thetastar_\alpha}{\sigma^2 I_\numobs}}
\Big[\inf_{\lambda > 0} \|\hat \theta_\lambda(X_\alpha, y) - \thetastar_\alpha\|_2^2\Big]
\geq 
\frac{c_2 c_1 p}{16}
\Big(
\frac{\sigma^2 B \dimension}{\numobs \alpha} 
\wedge 
R^2 \Big(\frac{\sigma^2 B}{R^2 \numobs} \alpha\Big)^{1 - p/2}
\wedge 
R^2
\Big)
\end{equation}
where above we have used $k \geq d/4$.

\paragraph{Lower bound when $p = 0$:} In this case, we
let $s' = \twomin{s}{k}$. Note that $s' \geq s/4$. We then set 
$\eta = 2 \sqrt{\tfrac{\sigma^2 B \alpha}{\numobs}} \1_{s'}$, 
and this yields the lower bound 
\[
T_2(\eta) \geq \frac{1}{8} \frac{\sigma^2 B s}{\numobs} 
\alpha \overline{W_{s'}}
\]
In this case, we have, after combining this bound with the bound on 
$T_1$ that for $\thetastar_\alpha \defn \thetastar(\eta)$ as defined above, 
\begin{equation}\label{ineq:final-lower-hard-sparse-lasso}
\E_{y \sim \Normal{X_\alpha \thetastar_\alpha}{\sigma^2 I_\numobs}}
\Big[\inf_{\lambda > 0} \|\hat \theta_\lambda(X_\alpha, y) - \thetastar_\alpha\|_2^2\Big]
\geq 
\frac{c_2 c_1 p}{16}
\Big(
\frac{\sigma^2 B \dimension}{\numobs \alpha} 
\wedge \frac{\sigma^2 B s}{\numobs} \alpha\Big)
\end{equation}

The proof of Theorem~\ref{thm:suboptimality-penalized-lasso} is 
complete after combining inequalities~\eqref{ineq:final-lower-weak-sparse-lasso}~\eqref{ineq:final-lower-hard-sparse-lasso}, and the following lemma.

\ble 
\label{lem:lasso-lower-constant-bound}
The constant factor $c \defn \tfrac{c_1 c_2 p}{16}$ 
is lower bounded as $c \geq \tfrac{9}{20000}$.
\ele 

We conclude this section by proving the lemmas above.

\subsubsection{Proof of Lemma~\ref{lem:lower-bound-one}}
\label{sec:proof-lem-lower-bound-one}
For the first term, $T_1$, we note that for any $\lambda \leq \lambdathresh$
we evidently have for each $i \in [k]$ and for any $\zeta > 0$, that
\[
\Big(\soft{\lambda B/\alpha}{z_i}\Big)^2 
= (|z_i| - \lambda B/\alpha)_+^2 
\geq (|z_i| - \lambdathresh B/\alpha)_+^2 
\geq 
\zeta^2 \frac{\sigma^2 B}{\numobs \alpha} 
\1\Big\{\tfrac{|z_i|}{\sqrt{\sigma^2 B/(n\alpha)}} \geq 1 + \zeta\Big\}
\]
Summing over $i \in [k]$, and taking $\zeta = 1/2$, we thus obtain the claimed almost sure lower bound.
\subsubsection{Proof of Lemma~\ref{lem:lower-bound-two}}
\label{sec:proof-lem-lower-bound-two}

Fix any $i$ such that $1 \leq i \leq k$. For any fixed
$\lambda \geq \lambdathresh$, note 
\[
S_{\lambda B}(z_{k + i}) \not \in [\eta_i/2, 3 \eta_i/2] 
\quad \mbox{implies} \quad 
\Big| S_{\lambda B}(z_{k + i}) - \eta_i \Big| \geq \frac{\eta_i}{2}. 
\]
Note that the condition 
$S_{\lambda B}(z_{k+i}) \not \in [\eta_i/2, 3 \eta_i/2]$ is 
equivalent to 
$z_{k + i} \not \in [\eta_i/2 + \lambda B, 3\eta_i/2 + \lambda B]$, 
Therefore, if $z_{k + i} \leq \eta_i/2 + \lambdathresh B$, then  
then for all $\lambda \geq \lambdathresh$ we have 
$|S_{\lambda B}(z_{k + i}) - \eta_i| \geq \frac{\eta_i}{2}$. 
Equivalently, we have that 
\begin{equation}\label{eqn:term-two-lower}
T_2 
\geq 
\frac{1}{4} \sum_{i=1}^k \eta_i^2 
\1\{z_{k+i} \leq \eta_i/2 + \lambdathresh B\} 
\geq 
\frac{1}{4} 
\sum_{i=1}^k \eta_i^2 \1\{z_{k + i} \leq \eta_i\}
\end{equation}
The final relation uses the distribution of $z_{k+i}$ and 
\[
\frac{-\eta_i/2 + \lambdathresh B}{\sqrt{\sigma^2 B/\numobs}}
= 
\sqrt{\alpha} - \frac{1}{2} \sqrt{\alpha} \sqrt{\frac{\numobs \eta_i^2}{\alpha \sigma^2 B}} 
\geq 0
\]
which holds by assumption that $\eta_i^2 \leq 4 \tfrac{\sigma^2 B}{n} \alpha$.

\subsubsection{Proof of Lemma~\ref{lem:lasso-lower-constant-bound}}

Evidently $c_1 \geq 1/2$ by symmetry. On the other hand, since $p \leq 1/2$, we have by anticoncentration results for 
Binomial random variables~\cite[Theorem 6.4]{GreMoh14} that 
$c_2 \geq p$. Therefore all together, $c \geq p^2/32$. Note that 
by standard lower bounds for the Gaussian tail~\cite[Theorem 1.2.6]{Dur19}, we have 
\[
p \geq \frac{10}{27} \e^{-9/8} \geq \frac{3}{25}, 
\]
which provides our claimed bound.

\subsection*{Acknowledgements}

RP gratefully acknowledges partial support from the ARCS Foundation via the Berkeley Fellowship. Additionally, RP thanks Martin J.\ Wainwright for 
helpful comments, discussion, and references in the preparation of this manuscript. CM  was partially supported by the National Science Foundation via grant DMS-2311127.

\bibliographystyle{plain}
\bibliography{references}

\begin{thebibliography}{10}

\bibitem{BirMas01}
Lucien Birg{\'e} and Pascal Massart.
\newblock Gaussian model selection.
\newblock {\em Journal of the European Mathematical Society}, 3:203--268, 2001.

\bibitem{buhlmann2011statistics}
Peter B{\"u}hlmann and Sara Van De~Geer.
\newblock {\em Statistics for high-dimensional data: methods, theory and
  applications}.
\newblock Springer Science \& Business Media, 2011.

\bibitem{candes2013well}
Emmanuel~J Candes and Mark~A Davenport.
\newblock How well can we estimate a sparse vector?
\newblock {\em Applied and Computational Harmonic Analysis}, 34(2):317--323,
  2013.

\bibitem{candes2009near}
Emmanuel~J. Cand\`es and Yaniv Plan.
\newblock Near-ideal model selection by {$\ell_1$} minimization.
\newblock {\em Ann. Statist.}, 37(5A):2145--2177, 2009.

\bibitem{DALALYAN17Lasso}
Arnak~S. Dalalyan, Mohamed Hebiri, and Johannes Lederer.
\newblock On the prediction performance of the {L}asso.
\newblock {\em Bernoulli}, 23(1):552--581, 2017.

\bibitem{Dic16}
Lee~H. Dicker.
\newblock Ridge regression and asymptotic minimax estimation over spheres of
  growing dimension.
\newblock {\em Bernoulli}, 22(1):1--37, 2016.

\bibitem{DonJoh94}
David~L. Donoho and Iain~M. Johnstone.
\newblock Minimax risk over {$\ell_p$}-balls for {$\ell_q$}-error.
\newblock {\em Probab. Theory Related Fields}, 99(2):277--303, 1994.

\bibitem{Dur19}
Rick Durrett.
\newblock {\em Probability---theory and examples}, volume~49 of {\em Cambridge
  Series in Statistical and Probabilistic Mathematics}.
\newblock Cambridge University Press, Cambridge, 2019.

\bibitem{eyre2023out}
Benjamin Eyre, Elliot Creager, David Madras, Vardan Papyan, and Richard Zemel.
\newblock Out of the ordinary: Spectrally adapting regression for covariate
  shift.
\newblock {\em arXiv preprint arXiv:2312.17463}, 2023.

\bibitem{fan2020statistical}
Jianqing Fan, Runze Li, Cun-Hui Zhang, and Hui Zou.
\newblock {\em Statistical foundations of data science}.
\newblock CRC press, 2020.

\bibitem{foygel2011fast}
Rina Foygel and Nathan Srebro.
\newblock Fast-rate and optimistic-rate error bounds for $\ell_1$-regularized
  regression.
\newblock {\em arXiv preprint arXiv:1108.0373}, 2011.

\bibitem{ge2023maximum}
Jiawei Ge, Shange Tang, Jianqing Fan, Cong Ma, and Chi Jin.
\newblock Maximum likelihood estimation is all you need for well-specified
  covariate shift.
\newblock {\em arXiv preprint arXiv:2311.15961}, 2023.

\bibitem{GreMoh14}
Spencer Greenberg and Mehryar Mohri.
\newblock Tight lower bound on the probability of a binomial exceeding its
  expectation.
\newblock {\em Statist. Probab. Lett.}, 86:91--98, 2014.

\bibitem{hastie2015statistical}
Trevor Hastie, Robert Tibshirani, and Martin Wainwright.
\newblock {\em Statistical learning with sparsity: the {L}asso and
  generalizations}.
\newblock CRC press, 2015.

\bibitem{hsu2012random}
Daniel Hsu, Sham~M Kakade, and Tong Zhang.
\newblock Random design analysis of ridge regression.
\newblock In {\em Conference on learning theory}, pages 9--1. JMLR Workshop and
  Conference Proceedings, 2012.

\bibitem{Johnstone2019}
Iain~M. Johnstone.
\newblock Gaussian estimation: Sequence and wavelet models.
\newblock Book manuscript, September 2019.

\bibitem{Kelner2021preconditioning}
Jonathan~A. Kelner, Frederic Koehler, Raghu Meka, and Dhruv Rohatgi.
\newblock On the power of preconditioning in sparse linear regression.
\newblock In {\em 2021 {IEEE} 62nd {A}nnual {S}ymposium on {F}oundations of
  {C}omputer {S}cience---{FOCS} 2021}, pages 550--561. 2022.

\bibitem{kpotufe2021marginal}
Samory Kpotufe and Guillaume Martinet.
\newblock Marginal singularity and the benefits of labels in covariate-shift.
\newblock {\em The Annals of Statistics}, 49(6):3299--3323, 2021.

\bibitem{lehmann2006theory}
Erich~L Lehmann and George Casella.
\newblock {\em Theory of point estimation}.
\newblock Springer Science \& Business Media, 2006.

\bibitem{lei2021near}
Qi~Lei, Wei Hu, and Jason Lee.
\newblock Near-optimal linear regression under distribution shift.
\newblock In {\em International Conference on Machine Learning}, pages
  6164--6174. PMLR, 2021.

\bibitem{ma2023optimally}
Cong Ma, Reese Pathak, and Martin~J Wainwright.
\newblock Optimally tackling covariate shift in {RKHS}-based nonparametric
  regression.
\newblock {\em The Annals of Statistics}, 51(2):738--761, 2023.

\bibitem{pathak2022new}
Reese Pathak, Cong Ma, and Martin Wainwright.
\newblock A new similarity measure for covariate shift with applications to
  nonparametric regression.
\newblock In {\em International Conference on Machine Learning}, pages
  17517--17530. PMLR, 2022.

\bibitem{Pat23}
Reese Pathak, Martin~J. Wainwright, and Lin Xiao.
\newblock Noisy recovery from random linear observations: Sharp minimax rates
  under elliptical constraints, 2023.

\bibitem{raskutti2011minimax}
Garvesh Raskutti, Martin~J Wainwright, and Bin Yu.
\newblock Minimax rates of estimation for high-dimensional linear regression
  over $\ell_q$-balls.
\newblock {\em IEEE transactions on information theory}, 57(10):6976--6994,
  2011.

\bibitem{shimodaira2000improving}
Hidetoshi Shimodaira.
\newblock Improving predictive inference under covariate shift by weighting the
  log-likelihood function.
\newblock {\em Journal of statistical planning and inference}, 90(2):227--244,
  2000.

\bibitem{tibshirani1996regression}
Robert Tibshirani.
\newblock Regression shrinkage and selection via the lasso.
\newblock {\em Journal of the Royal Statistical Society Series B: Statistical
  Methodology}, 58(1):267--288, 1996.

\bibitem{van2018tight}
Sara van~de Geer.
\newblock On tight bounds for the {L}asso.
\newblock {\em J. Mach. Learn. Res.}, 19:Paper No. 46, 48, 2018.

\bibitem{wainwright2019high}
Martin~J Wainwright.
\newblock {\em High-dimensional statistics: A non-asymptotic viewpoint},
  volume~48.
\newblock Cambridge university press, 2019.

\bibitem{wang2023pseudo}
Kaizheng Wang.
\newblock Pseudo-labeling for kernel ridge regression under covariate shift.
\newblock {\em arXiv preprint arXiv:2302.10160}, 2023.

\bibitem{zhang2022class}
Xuhui Zhang, Jose Blanchet, Soumyadip Ghosh, and Mark~S Squillante.
\newblock A class of geometric structures in transfer learning: Minimax bounds
  and optimality.
\newblock In {\em International Conference on Artificial Intelligence and
  Statistics}, pages 3794--3820. PMLR, 2022.

\end{thebibliography}
\appendix

\section{Results in the Gaussian sequence model}
\label{app:gaussian-sequence-model}
In this section, we collect classical results regarding the nonasymptotic 
minimax rate of estimation for Gaussian sequence model over 
the origin-centered $\ell_p$ balls, $p \in [0, 1]$. 
All of the results in this section are based on the monograph~\cite{Johnstone2019}.
We use the following notation to specify the minimax rate of interest, 
\[
\mathfrak{M}\Big(p, d, R, \eps\Big)
\defn 
\inf_{\hat \mu}
\sup_{\substack{\mu^\star \in \R^\dimension \\ 
\|\mu^\star\|_p \leq R}}
\E_{y \sim \Normal{\mu^\star}{\eps^2}} 
\Big[\twonorm{\hat \mu(y) - \mu^\star}^2\Big].
\]
As usual, the infimum ranges over measurable estimators from 
the observation vector $y \in \R^\dimension$ to an estimate $\hat \mu(y) 
\in \R^\dimension$. Throughout, we use the notation $\tau \defn \frac{\eps}{R}$
for the inverse signal-to-noise ratio. 

\bpr [Minimax rate of estimation when $0 < p \leq 1$]
\label{prop:minimax-rate-weak-sparse-gsm}
Fix an integer $d \geq 1$. 
Let $p \in (0, 1]$. If $R, \eps > 0$ satisfy 
\[
\frac{1}{d^{2/p}} \leq \tau^2 \leq \frac{1}{1 + \log \dimension}
,
\qquad \mbox{where} \quad  \tau = \frac{\eps}{R},
\]
then 
\[
\frac{7}{2000}\, R^2 \, (\invsnr^2 \log(\e d \invsnr^p))^{1 - p/2}
\leq 
\mathfrak{M}(p, d, R, \eps)
\leq 
1203 \,  R^2 \, (\invsnr^2 \log(\e d \invsnr^p))^{1 - p/2}.
\]
\epr 

The upper and lower bounds are taken from Theorem 11.7 in 
the monograph~\cite{Johnstone2019}. 
Although the constants are not made explicit in their theorem statement, 
the upper bound constant is obtained via their Theorem 11.4, 
setting their parameters as $\zeta = \tfrac{23}{4}, \gamma = 2\e, \beta = 0$. 
Similarly, the lower bound constant is implicit in their proof of Theorem 11.7. 

We now turn to the minimax rate in the special case that $p = 0$. 

\bpr [Minimax rate of estimation when $p = 0$]
\label{prop:minimax-rate-hard-sparse-gsm}
Suppose that $d \geq 1$ and $s \in [d]$. 
Then for any $\eps > 0$ we have 
\[
\frac{3}{500} \, 
\eps^2 \, s \log \Big(\e \frac{d}{s}\Big)
\leq 
\mathfrak{M}(p, d, s, \eps)
\leq 
2 \,
\eps^2 \, s \log \Big(\e \frac{d}{s}\Big),
\]
provided that $p = 0$.
\epr 

The proof of the above claim is omitted as it is a straightforward combination of the standard minimax rate $\eps^2 k$ for 
the unconstrained Normal location model in a $k$-dimensional problem (this provides a useful lower bound when $s \geq d/2$ or when $d = 1$) and the result in Proposition 8.20 in the monograph~\cite{Johnstone2019}.

\section{Details for experiments in Figure~\ref{fig:simulation}}
\label{app:sim-details}

For each choice of $p$, we simulate the oracle Lasso and
STOLS procedures on instances $(X_\numobs, \theta_\numobs^\star)$ indexed by the sample
size $n \in \{1000, 2000, 3000, 5000, 10000, 15000\}$.
The matrix $X_\numobs \in \R^{\numobs \times \numobs}$ is
block diagonal and given by
\[
\frac{1}{\sqrt{\numobs}} X_\numobs = \begin{pmatrix}
  I_{n/2 \times n/2} & 0 \\
  0 & n^{-1/4} I_{n/2 \times n/2} 
\end{pmatrix},
\]
When $p = 0$, $\thetastar_\numobs = 2 e_{n/2 +1}$ and when
$p \neq 0$, $\thetastar_\numobs = e_{n/2 +1}$. 
In the figures, we
are plotting the average performance of the oracle Lasso and STOLS procedures,
as measured by $\ell_2$ error, when applied to the data $(X_\numobs, y)$,
where $y \sim \Normal{X_\numobs \thetastar_\numobs}{I_\numobs}$.
The average is taken over $1000$ trials for $n < 10,000$. In the case 
$\numobs \geq 10,000$ due to memory constraints we only run $300$ trials.

The STOLS procedure
is implemented as described in Section~\ref{sec:STOLS}.
On the other hand, the oracle Lasso procedure is implemented by a slightly
more involved procedure. Our goal is to compute
\[
\hat \theta_{\hat \lambda}(X, y)
\quad \mbox{where} \quad
\hat \lambda \in \argmin_{\lambda > 0} \|\hat \theta_\lambda(X, y) - \thetastar\|_2^2,
\]
where the Lasso is defined as in display~\eqref{eqn:def-lasso}. To do this,
we can use the fact that the Lasso regularization path is piecewise linear.
That is, there exist knot points $0 = \lambda_0 < \lambda_1 < \lambda_2 < \dots < \lambda_m$ such that
the knot points $\hat \theta_i \defn \hat \theta_{\lambda_i}(X, y)$ satisfy
$\|\theta_i\|_0 > \|\theta_{i+1}\|_0$. Moreover,
we have
\[
\{ \hat \theta_{\lambda}(X, y) : \lambda \in (\lambda_{i}, \lambda_{i+1}) \}
=
\{ \hat \theta_i + \alpha (\hat \theta_{i+1} - \hat \theta_i) : \alpha \in (0, 1) \}. 
\]
That is, we can compute the set of Lasso solutions between the knot points by taking
all convex combinations of knot points. Therefore the distance between the
oracle Lasso solution and the true parameter $\thetastar$ satisfies,
\[
\|\hat \theta_{\hat \lambda}(X, y) - \thetastar\|_2^2
= \min_{i} \min_{\alpha \in [0, 1]}
\|\hat \theta_i + \alpha(\hat \theta_{i+1} - \hat \theta_i) - \thetastar\|_2^2.
\]
We are able to compute the righthand side of the display above by noting that for
each $i$ the inner minimization problem is a quadratic function of the univariate parameter
$\alpha$ and therefore can be minimized explicitly.

\paragraph{Code:} The code has been released at the following public repository,
\begin{quote}
  \texttt{https://github.com/reesepathak/lowerlassosim}.
\end{quote}
In particular, the repository contains a \texttt{Python} program which runs
simulations of STOLS and oracle Lasso on the lower bound instance described above
for any desired choice of $p \in [0, 1]$. 

\end{document}